\documentclass{amsart}

\usepackage{graphicx}
\usepackage{amsfonts, amsmath, amscd, a4wide}
\usepackage{a4wide}
\usepackage{psfrag}
\usepackage[all]{xy}

\let\cal\mathcal
\def\AA{{\cal A}}
\def\BB{{\cal B}}
\def\CC{{\cal C}}

\def\EE{{\cal E}}
\def\FF{{\cal F}}
\def\GG{{\cal G}}
\def\HH{{\cal H}}

\def\OO{{\cal O}}

\def\TT{{\cal T}}

\let\blb\mathbb

\def\bQ{{\blb Q}}

\def\bZ{{\blb Z}}

\def\bN{{\blb N}}
\def\bR{{\blb R}}

\def\bZ{{\blb Z}}

\def\Mod{\operatorname{Mod}}
\def\Modfd{\operatorname{\Mod^{\text{fd}}}}

\def\qgr{\operatorname{qgr}}

\def\coh{\mathop{\text{\upshape{coh}}}}

\def\rad{\operatorname {rad}}

\def\Ext{\operatorname {Ext}}

\def\Hom{\operatorname {Hom}}

\def\End{\operatorname {End}}
\def\RHom{\operatorname {RHom}}

\def\CYdim{\operatorname {CYdim}}
\def\GKdim{\operatorname {GKdim}}
\def\gldim{\operatorname {gl.dim}}

\def\im{\operatorname {im}}

\def\cone{\operatorname {cone}}
\def\coker{\operatorname {coker}}
\def\ker{\operatorname {ker}}
\def\Ker{\operatorname {ker}}

\def\End{\operatorname {End}}

\def\rk{\operatorname {rk}}

\def\gldim{\operatorname {gl\,dim}}

\def\exist{\exists}

\def\dom{\operatorname {dom}}
\def\cohproj{\operatorname {cohproj}}

\DeclareMathOperator{\ind}{ind}

\newcommand\Db{\operatorname{D^{b}}}

\newtheorem{lemma}{Lemma}[section]
\newtheorem{proposition}[lemma]{Proposition}
\newtheorem{theorem}[lemma]{Theorem}
\newtheorem{corollary}[lemma]{Corollary}

\theoremstyle{definition}

\newtheorem{example}[lemma]{Example}
\newtheorem{definition}[lemma]{Definition}

{

}

\theoremstyle{remark}

\newtheorem{remark}[lemma]{Remark}

\newdimen\uboxsep \uboxsep=1ex
\def\uboxn#1{\vtop to 0pt{\hrule height 0pt depth 0pt\vskip\uboxsep
\hbox to 0pt{\hss #1\hss}\vss}}

\def\uboxs#1{\vbox to 0pt{\vss\hbox to 0pt{\hss #1\hss}
\vskip\uboxsep\hrule height 0pt depth 0pt}}

\def\Ob{\operatorname{Ob}}

\newcommand\exa{\nopagebreak \begin{center}\smallskip \nopagebreak               \begin{minipage}[t]{6cm}\sloppy}
\newcommand\exb{\end{minipage}\kern 1cm\begin{minipage}[t]{8cm}\sloppy}
\newcommand\exc{\end{minipage}\kern -3cm \smallskip\end{center}}

\title{Abelian 1-Calabi-Yau Categories}
\author{Adam-Christiaan van Roosmalen}
\address{Adam-Christiaan van Roosmalen\\Hasselt University
\\Research group Algebra\\Agoralaan, gebouw D\\B-3590 Diepenbeek (Belgium)}\email{AdamChristiaan.vanRoosmalen@UHasselt.be}

\begin{document}

\begin{abstract}
In this paper, we show all $k$-linear abelian 1-Calabi-Yau categories over an algebraically closed field $k$ are derived equivalent to either the category of coherent sheaves on an elliptic curve, or to the finite dimensional representations of $k[[t]]$.  Since all abelian categories derived equivalent with these two are known, we obtain a classification of all $k$-linear abelian 1-Calabi-Yau categories up to equivalence.
\end{abstract}

\maketitle

\section{Introduction}

In this paper, we classify  \emph{abelian 1-Calabi-Yau categories}
over an algebraically closed field~$k$.  Recall that an abelian
1-Calabi-Yau category is a $k$-linear $\Hom/\Ext$-finite abelian category
together with
natural isomorphisms $\Hom(X,Y) \cong \Ext(Y,X)^\ast$ for $X,Y\in \AA$.

Our main result (reformulated in the body of the text as Theorem
\ref{theorem:Main}) is the following.
\begin{theorem}
\label{mainth}
Let $\mathcal{A}$ be an indecomposable abelian 1-Calabi-Yau
   category. Then $\mathcal{A}$ is \emph{derived} equivalent to one of the
 following two categories.
\begin{enumerate}
\item  Finite dimensional
representations of $k[[t]]$.
\item The category of coherent sheaves on an elliptic
curve.
\end{enumerate}
\end{theorem}
There is a general interest in the classification of categories which
are homologically small in some sense (see e.g.\ \cite{Happel01},  \cite{Keller06}
\cite{ReVdB02}, \cite{vanRoosmalen06}). The above
theorem represents an enhancement of our knowledge in this area.

Besides this general motivation we mention the following particular
application. Recently Polishchuk and Schwartz \cite{Polishchuk03} constructed a
category $\CC$ of holomorphic vector bundles on a
\emph{non-commutative 2-torus}.  Polishchuk subsequently showed that
$\CC$ is derived equivalent to the category of coherent sheaves on an
elliptic curve \cite{Polishchuk03}.  Part of Polishchuk's proof amounts to
establishing the highly non-trivial fact that $\CC$ is $1$-Calabi-Yau
\cite[Cor 2.12]{Polishchuk03}. Once one knows this, one could now finish the
proof by simply invoking Theorem \ref{mainth} (with $\AA$ being a
suitable abelian hull of~$\CC$).

\medskip

We  briefly outline some steps in the proof of Theorem \ref{mainth}.  Some of
our tools come from  representation theory of algebras and non-commutative
algebraic geometry.  Other
tools were already employed by Polishchuk, but are now used in a more
abstract setting.

Fix a connected abelian 1-Calabi-Yau category $\AA$.  First, we prove
the existence of \emph{endo-simple} objects in $\AA$, i.e.\ objects
$X\in \AA$ such that $\End X \cong k$. Associated to such objects
there are \emph{twist functors} \cite{Seidel01} $T_A,
T_A^\ast$. These functors are mutually
inverse auto-equivalences of $D^b(\AA)$ which on objects take the values
$T_X Y =
\cone (X \otimes \RHom(X,Y) \longrightarrow Y)$ and $T_X^\ast Y =
\cone (Y[-1] \longrightarrow X[-1] \otimes \RHom(Y,X)^\ast)$.

Using twist functors we establish various useful facts. Most notably, we
prove that the subcategory of endo-simple objects in $\AA$ has no cycles of
non-zero maps (Proposition \ref{proposition:EndoSimplesDirected}) and
hence is ordered.  We also show that all Auslander-Reiten components
of $\AA$ are homogeneous tubes based on endo-simple objects (Proposition
\ref{proposition:Building}).

We may assume that $\AA$ has at least two non-isomorphic endo-simple
objects as the remaining case is easily disposed with. Using
connectedness and the results mentioned in the previous paragraph we
may in fact select non-isomorphic endo-simple objects $E$ and $B$ such
that $\Hom(E,B)\neq 0$. After doing so we consider the sequence of
objects $\EE=(T_B^{n} E)_{n \in \bZ}$ in $\Db \AA$. We construct a
certain associated $t$-structure on $D^b(\AA)$ with heart $\HH$ such
that $\EE$ is an ample sequence in the sense of \cite{Polishchuk05} in
$\HH$. Hence $\EE$ defines a finitely presented graded coherent
algebra $A$ such that $\HH$ is equivalent to the category
$\operatorname{qgr}(A)$ of finitely presented graded $A$-modules modulo
finite dimensional ones.

We then show that $A$ is a domain of Gelfand-Kirillov dimension two
and we invoke the celebrated Artin and Stafford classification theorem
\cite{Artin95} which shows that $\operatorname{qgr}(A)$ is of the form
$\operatorname{coh}(X)$ for a projective curve $X$. Since $\HH$ is
1-Calabi-Yau this implies that $X$ must be an elliptic curve,
finishing the proof.

It is not hard to describe the abelian 1-Calabi-Yau
categories that occur within the derived equivalence classes in
Theorem \ref{mainth} (see e.g. \cite{Burban06}).  We discuss this using the
language of this paper in \S\ref{subsection:Abelian}.

\subsection*{Acknowledgment}
The author thanks Michel Van den Bergh for useful discussions as well as
for contributing some ideas.
\section{Preliminaries}

Throughout this paper, fix an algebraically closed field $k$ of arbitrary characteristic.  All algebras and categories are assumed to be $k$-linear.

We will also assume all abelian categories are \emph{connected} in the sense that between two indecomposable objects there is an unoriented path of non-zero maps between indecomposables.

If $\AA$ is abelian, we write $\Db \AA$ for the bounded derived category of $\AA$.  The category $\Db \AA$ has the structure of a triangulated category.  Whenever we use the word "triangle" we mean "distinguished triangle".

An abelian or triangulated category $\AA$ is \emph{Ext-finite} if for all objects $X,Y \in \Ob(\AA)$ one has that $\dim_{k}\Ext^{i}(X,Y)<\infty$ for all $i \in \bN$.  We say that $\AA$ is \emph{hereditary} if $\Ext^{i}(X,Y)=0$ for all $i \geq 2$.

\subsection{Serre duality}
If $\CC$ is triangulated category, we will say that $\CC$ satisfies \emph{Serre duality} if there exists an auto-equivalence $F:\CC \to \CC$, called the \emph{Serre functor}, such that, for all $X,Y\in\Ob \CC$, there is an isomorphism
$$\Hom_{\CC}(X,Y) \cong \Hom_{\CC}(Y,FX)^{*}$$
which is natural in $X$ and $Y$ and where $(-)^{*}$ denotes the vector space dual.

We will say an abelian category $\AA$ has Serre duality if the category $\Db \AA$ has a Serre functor.

It has been proven in \cite{ReVdB02} that an abelian category $\AA$ without non-zero projectives has a Serre functor if and only if the category $\Db\AA$ has Auslander-Reiten triangles.  In this case the action of the Serre functor on objects coincides with $\tau[1]$, where $\tau$ is the Auslander-Reiten translation.

\subsection{Calabi-Yau categories}
Let $\AA$ be an Ext-finite abelian category with Serre duality.  We will say $\AA$ is \emph{Calabi-Yau of dimension $n$} or shorter that $\AA$ is \emph{$n$-Calabi-Yau} if $F \cong [n]$ for a certain $n \in \bN$, thus if the $n^\text{th}$ shift is a Serre functor.  We write $\CYdim \AA = n$.

The following well-known property relates the Calabi-Yau dimension and the homological dimension.

\begin{proposition}\label{proposition:Dimension}
Let $\AA$ be an abelian Calabi-Yau category.  Then $\CYdim \AA = \gldim \AA$.
\end{proposition}

\begin{proof}
Let $n = \CYdim \AA$, then for every $X \in \Ob \AA$ we have $\Hom(X,X) \cong \Ext^n(X,X)^{*}$.  Since the former is non-zero, we see $\CYdim \AA \leq \gldim \AA$.

Let $i \in \bN$ and $X,Y \in \Ob \AA$ be chosen such that $\Ext^{i}(X,Y) \not = 0$.  Using the Calabi-Yau property, we find $\Ext^{i}(X,Y) \cong \Ext^{n-i}(Y,X)$, hence $n \geq i$.  We find $\CYdim \AA \geq \gldim \AA$.
\end{proof}

In particular, if $\AA$ is a 1-Calabi-Yau category, then $\AA$ is hereditary.  Since $F \cong [1]$ and $F$ coincides with $\tau[1]$ on indecomposables of $\AA$, it follows that $\tau$ is naturally isomorphic to the identity functor on $\AA$ and hence $\Hom_\AA(X,Y) \cong \Ext_\AA(Y,X)^\ast$, for all objects $X,Y \in \AA$.

\subsection{Twist functors}
Let $\AA$ be an abelian 1-Calabi-Yau category.  For an object $A \in \Db \AA$, we may consider the \emph{twist functors}, $T_A$ and $T_A^\ast$, in $\Db \AA$ whose values on objects are up to isomorphism characterized by the following triangles
$$T_A X[-1] \longrightarrow A \otimes \RHom(A,X) \stackrel{\epsilon}{\longrightarrow} X \longrightarrow T_A X$$
and
$$T_A^\ast X \longrightarrow X \stackrel{\epsilon^\ast}{\longrightarrow} A \otimes \RHom(X,A)^\ast \longrightarrow T_A^\ast X[1]$$
where $\epsilon : A \otimes \RHom(A,X) \longrightarrow X$ and $\epsilon^\ast : X \longrightarrow A \otimes \RHom(X,A)^\ast$ are the canonical morphisms.

Let $S$ be an \emph{endo-simple} object, i.e.\ $\End(S) \cong k$.  Since $\AA$ is 1-Calabi-Yau, we know from \cite[Proposition 2.10]{Seidel01} that $T_S$ and $T_S^\ast$ are inverses.  In particular, they are autoequivalences.

\subsection{Ample sequences}
For the benefit of the reader, we will recall some definitions and results from \cite{Polishchuk05} which will be used in the rest of this paper.  Throughout, let $\AA$ be a Hom-finite abelian category.

We begin with the definition of ample sequences.
\begin{enumerate}
\item A sequence $\EE = (E_i)_{i\in \bZ}$ is called \emph{projective} if for every epimorphism $X \to Y$ in $\AA$ there is an $n\in \bZ$ such that $\Hom(E_i,X) \to \Hom(E_i,Y)$ is surjective for $i < n$. 

\item A projective sequence $\EE = (E_i)_{i\in \bZ}$ is called \emph{coherent} if for every $X \in \Ob \AA$ and $n \in \bZ$, there are integers $i_1, \ldots, i_s \leq n$ such that the canonical map
$$\bigoplus_{j=1}^s \Hom(E_i, E_{i_j}) \otimes \Hom(E_{i_j},X) \longrightarrow \Hom(E_i,X)$$
is surjective for $i<<0$.

\item A coherent sequence $\EE =  (E_i)_{i\in \bZ}$ is \emph{ample} if for all $X \in \AA$ the map $\Hom(E_i,X) \not=0$ for $i<<0$.
\end{enumerate}

Let $A_{ij} = \Hom(E_i,E_j)$ for $i \leq j$.  We may define an algebra $A = A(\EE) = \oplus_{i \leq j} A_{ij}$ in a natural way.  If $A_{ii} \cong k$, then $A$ is a \emph{coherent $\bZ$-algebra} in the sense of \cite{Polishchuk05} (see \cite[Proposition 2.3]{Polishchuk05}).

We will refer to the right $A$-modules having a resolution by finitely generated projectives as \emph{coherent modules}.  These modules form an abelian category, $\coh A$, and the finite dimensional modules form a Serre subcategory denoted by $\coh^{b} A$.  We define the quotient 
$$\cohproj A \cong \coh A / {\coh}^{b} A.$$

We may use this to give a description of Ext-finite abelian categories with an ample sequence.

\begin{theorem}\label{theorem:Polishchuk}
\cite[Theorem 2.4]{Polishchuk05} Let $\EE=(E_i)$ be an ample sequence, $A=A(\EE)$ the corresponding $\bZ$-algebra, then there is a equivalence of categories $\AA \cong \cohproj A$.
\end{theorem}

We will be interested in the special case where there is an automorphism $t : \Db \AA \longrightarrow \Db \AA$ such that $E_i \cong t^i E$.  We let $R = R(\EE) = \oplus_{i \in \bN} R_i$ where $R_i = \Hom(E,t^i E)$ and make it into a $\bZ$-graded algebra in an obvious way.

If $R$ is noetherian then the coherent $R$-modules correspond to the finitely generated ones and $\cohproj R$ corresponds to $\qgr R$, the finitely generated modules modulo the finite dimensional ones.

We will use following corollary of Theorem \ref{theorem:Polishchuk}.

\begin{corollary}\label{corollary:Polishchuk}
Let $\AA$ be a Hom-finite abelian category, $t$ be an autoequivalence of $\AA$ and $E$ an object of $\AA$.  If $\EE=(t^i E)$ is an ample sequence and the corresponding graded algebra $R = R(\EE)$ is noetherian, then $\AA \cong \qgr R$.
\end{corollary}

\subsection{$t$-structures}

In order to find derived equivalent categories, we will use the theory of $t$-structures \cite{Beilinson82}.

\begin{definition}\label{definition:t}
A \emph{$t$-structure} on a triangulated category $\CC$ is a pair $(D^{\geq 0}, D^{\leq 0})$ of non-zero full subcategories of $\CC$ satisfying the following conditions, where we denote $D^{\leq n} = D^{\leq 0} [-n]$ and $D^{\geq n} = D^{\geq 0} [-n]$
\begin{enumerate}
\item $D^{\leq 0} \subseteq D^{\leq 1}$ and $D^{\geq 1} \subseteq D^{\geq 0}$
\item $\Hom(D^{\leq 0}, D^{\geq 1}) = 0$
\item\label{split} $\forall Y \in \CC$, there exists a triangle $X \to Y \to Z \to X[1]$ with $X \in D^{\leq 0}$ and $Z \in D^{\geq 1}$.
\end{enumerate}
\end{definition}
Furthermore, we will say the $t$-structure is \emph{bounded} if $\bigcap_n D^{\leq n} = \bigcap_n D^{\geq n} = \{0\}$.

We will say a $t$-structure is \emph{split} if all triangles in (\ref{split}) are split, or equivalently, if $\ind \CC = \ind D^{\geq 1} \cup \ind D^{\leq 0}$.  We have following result.

\begin{theorem}\cite{Berg}\label{theorem:Berg}
Let $\AA$ be an abelian category and let $(D^{\geq 0}, D^{\leq 0})$ be a bounded $t$-structure on $\Db \AA$.  Then the heart $\HH$ is hereditary if and only if $(D^{\geq 0}, D^{\leq 0})$ is a split $t$-structure.  In this case, $\AA$ and $\HH$ are derived equivalent.
\end{theorem}

\subsection{Elliptic curves}

For the benefit of the reader, we recall certain properties of the category of coherent sheaves on an elliptic curve $X$.  This category has first been described in \cite{Atiyah57}.

An elliptic curve is a curve of genus 1 and thus, in particular, $\AA = \coh X$ is a 1-Calabi-Yau category.

Let $\OO$ be the structure sheaf and, for a point $P$, let $k(P)$ be a torsion sheaf.  For a coherent sheaf $\FF$ the \emph{degree} and \emph{rank} may be defined as 
\begin{eqnarray*}
\deg \EE &=& \chi (\OO,\EE) \\
\rk \EE &=& \chi(\EE,k(P)),
\end{eqnarray*}
respectively.  It follows from the Riemann-Roch theorem that
\begin{equation}\label{equation:RR}
\chi(\EE,\FF) = \deg \FF \rk \EE - \deg \EE \rk \FF.\end{equation}

Furthermore, the \emph{slope} of $\EE$ is defined as $\mu(\EE)=\frac{\deg \EE}{\rk \EE} \in \bQ\cup\{ \infty \}$.  A coherent sheaf $\FF$ is called \emph{stable} or \emph{semi-stable} if for every short exact sequence $0 \longrightarrow \EE \longrightarrow \FF \longrightarrow \GG \longrightarrow 0$ we have $\mu(\EE) \leq \mu(\FF)$ or $\mu(\EE) < \mu(\FF)$, respectively.

It is well-known that all indecomposable coherent sheaves are semi-stable.  For stable sheaves, we have the following equivalent conditions
\begin{enumerate}
\item $\EE$ is stable,
\item $\EE$ is endo-simple, i.e. $\End(\EE) \cong k$,
\item $\rk \EE$ and $\deg \EE$ are coprime.
\end{enumerate}

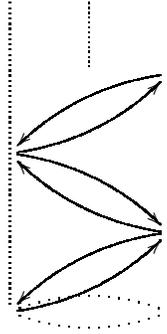
\begin{figure}[tbp]
	\centering

		$\xymatrix{& & \\
				& \ar@{.}[u] & \ar@/_/[lld] \\
				\ar@/_/[rru] \ar@/^/[rrd] && \\
				&& \ar@/_/[lld] \ar@/^/[llu] \\
		  	\ar@/_/[rru] \ar@{.}@/_/[rr] \ar@{.}@/^/[rr] \ar@{.}[uuuu] & & \ar@{.}[uuuu]}$
		
	\caption{A homogeneous tube.}
	\label{figure:Tube}
\end{figure}

Every semi-stable sheaf is a finite extension of an endo-simple one with itself.  We may visualise this via the Auslander-Reiten quiver of $\coh X$.  All Auslander-Reiten components are homogeneous tubes, i.e. components of the form $\bZ A_\infty / \langle \tau \rangle$, cfr. Figure \ref{figure:Tube}, where the bottom element is a stable sheaf.

Every such tube corresponds to an abelian subcategory of $\coh X$ equivalent to $\Modfd k[[t]]$ and all indecomposable objects in the same homogeneous tube have the same slope.  Thus the full subcategory of $\coh X$ spanned by all indecomposable objects of a given slope is an abelian subcategory of $\coh X$ and is of the form $\oplus \Modfd k[[t]]$, where the sum is indexed over the stable objects with the given slope.

Finally, it follows directly from (\ref{equation:RR}) that, for non-isomorphic stable sheaves, $\EE$ and $\FF$, we have $\Hom(\EE,\FF) \not= 0$ if and only if $\mu(\EE) < \mu(\FF)$.  Thus for semi-stable sheaves $\EE'$ and $\FF'$ we have $\Hom(\EE',\FF') \not= 0$ if and only if $\mu(\EE') < \mu(\FF')$ or $\EE'$ and $\FF'$ lie in the same tube.
\section{Endo-simple objects}

Let $\AA$ be a connected $k$-linear abelian 1-Calabi-Yau category.  It will turn out that the endo-simple objects are the building blocks of $\AA$.  Therefore, in this section, we will give some properties of endo-simple objects.  Recall that $X$ is an endo-simple object if $\End X \cong k$.  It follows from the Calabi-Yau property that every endo-simple object is 1-spherical in the sense of \cite{Seidel01}.

\begin{proposition}\label{proposition:EndoSimple}
Let $\CC$ be a Hom-finite abelian category.  For every object $X \in \Ob\CC$ there exists an endo-simple object occurring both as subobject and quotient object of $X$.  In particular, $\CC$ has an endo-simple object.
\end{proposition}

\begin{proof}
Assume $X$ is not endo-simple and let $f : X \to X$ be a non-invertible endomorphism.  We show that $\dim \End I < \dim \End X$ where $I = \im f$.

Indeed, since we have an epimorphism $X \to I$ and monomorphism $I \to X$, we get a composition of monomorphisms $\Hom(I,I) \to \Hom(X,I) \to \Hom(X,X)$.  Since the image of this composition has to be in $\rad(X,X)$, we know $\dim \Hom(I,I) < \dim \Hom(X,X)$.  Iteration finishes the proof.
\end{proof}

\begin{proposition}\label{proposition:CanonicalMap}
Let $S$ be an endo-simple object and $X \in \ind \AA$.  Each of the canonical maps $S \otimes \Hom(S,X) \longrightarrow X$ and $X \longrightarrow S \otimes \Hom(X,S)^\ast$ is either a monomorphism or an epimorphism.  If $\Hom(X,S) \not= 0$, then the first map is a monomorphism.  If $\Hom(S,X) \not= 0$, then the latter is an epimorphism.
\end{proposition}

\begin{proof}
Consider in the derived category $\Db \AA$ the twist functor $T_S$ characterized by
$$T_S X[-1] \longrightarrow S \otimes \RHom(S,X) \stackrel{\epsilon}{\longrightarrow} X \longrightarrow T_S X.$$
It is shown in \cite{Seidel01} that this is an equivalence.  Applying the homological functor $H^0$ gives the long exact sequence
$$0 \to H^{-1}(T_S X) \to S \otimes \Hom(S,X) \stackrel{H^0 \epsilon}{\longrightarrow} X \to H^0(T_S X) \to S \otimes \Ext(S,X) \to 0.$$
Since $X$ is indecomposable and $T_S$ is an equivalence, either $H^{-1}(T_S X)$ or $H^0(T_S X)$ is zero, hence $H^0 \epsilon$ is a monomorphism or an epimorphism, respectively.

If we assume furthermore $\Hom(X,S) \not= 0$, and hence by the Calabi-Yau property $\Ext(S,X) \not= 0$, we find $H^0(T_S X) \not= 0$.  Hence $H^{-1}(T_S X) = 0$ and the canonical map $S \otimes \Hom(S,X) \to X$ is a monomorphism.

The other case is dual.
\end{proof}

\begin{proposition}\label{proposition:EndoSimplesDirected}
The subcategory of endo-simples is a directed category.
\end{proposition}

\begin{proof}
Let $S_0 \to S_1 \to \cdots \to S_n \to S_0$ be a cycle of non-zero non-isomorphisms between endo-simple objects.  We will assume $n$ is minimal with the property that such a cycle exists.

By Proposition \ref{proposition:CanonicalMap} we know the canonical map $\epsilon:S_0 \otimes \Hom(S_0,S_1) \to S_1$ is either a monomorphism or an epimorphism.  If $\epsilon$ is a monomorphism, then we know the composition
$$S_n \otimes \Hom(S_0,S_1) \longrightarrow S_0 \otimes \Hom(S_0,S_1) \stackrel{\epsilon}{\longrightarrow} S_1$$
is non-zero.  This induces a non-zero morphism $f: S_n \longrightarrow S_1$.  Since $f$ factors through $S_0 \otimes \Hom(S_0,S_1)$, we know $f$ is not invertible. 

Likewise, if $\epsilon$ is an epimorphism, we find a non-zero non-invertible morphism $S_0 \to S_2$.  In both cases we have found a shorter cycle, contradicting with the minimality of $n$.
\end{proof}

We now wish to show that every object has a composition series with endo-simple quotients.  Even more so, every indecomposable object has a composition series in which only one isomorphism class of endo-simple objects occur.  We start with a lemma.

\begin{lemma}\label{lemma:Cone}
Let $X \in \ind \AA$ such that the endo-simple object $S$ occurs both as subobject and quotient object of $X$.  If $C = \coker(S \otimes\Hom(S,X) \longrightarrow X)$ is not zero, then S occurs as both subobject and quotient object of every direct summand of $C$.
\end{lemma}

\begin{proof}
Assume $C \not\cong 0$.  Consider the exact sequence
$$0 \to S \otimes \Hom(S,X) \longrightarrow X \longrightarrow H^0(T_S X) \longrightarrow S \otimes \Ext(S,X) \to 0.$$
from the proof of Proposition \ref{proposition:CanonicalMap}. We may splice this as
$$0 \longrightarrow S \otimes \Hom(S,X) \longrightarrow X \longrightarrow C \longrightarrow 0$$
and
$$0 \longrightarrow C \longrightarrow H^0(T_S X) \longrightarrow S \otimes \Ext(S,X) \longrightarrow 0.$$

Since $T_S$ is an automorphism and $X$ is indecomposable, we know $H^0(T_S X)$ is indecomposable.  It now follows directly from \cite[Lemma 2*]{Ringel05} that $\Hom(S,C_1) \cong \Ext(C_1,S)^* \not=0$ and $\Hom(C_1,S) \cong \Ext(S,C_1)^* \not=0$ for every direct summand $C_1$ of $C$.  Proposition \ref{proposition:CanonicalMap} now yields that $S$ is both a subobject and a quotient object of every direct summand of $C$.
\end{proof}

\begin{proposition}\label{proposition:Building}
Every indecomposable object is obtained by repeatedly extending a given endo-simple with itself.
\end{proposition}

\begin{proof}
Let $S$ be an endo-simple object and denote by $\AA_S$ the full subcategory of $\AA$ spanned by the objects $Z$ which can be obtained from $S$ by taking a finite amount of extensions with itself.  The number of such extensions needed, will be denoted by $l_S(Z)$, and we will refer to it as the length of $Z$.

Since $\AA_S$ is a hereditary category with a unique simple $S$ such that $\dim \Ext(S,S)=1$, it follows easily that $\AA_S$ is equivalent to the finite dimensional representations of $k[[t]]$.

We will prove that if $X$ is an indecomposable object of $\AA$ such that $S$ occurs both as quotient and subobject, then $X \in \AA_S$.  Note that by Proposition \ref{proposition:EndoSimple} we may assume such an $S$ exists.

For every subobject $A$ of $X$ in $\AA_S$ and quotient object $B$ of $X$ in $\AA_S$, we have
\begin{equation*}
\dim \End_\AA X \ge \min (l_S(A),l_S(B)),
\end{equation*}
thus we may deduce either the length of such subobjects or the length of such quotient objects is bounded.  Assume that the length of $A$ is bounded, the other case is dual.

We will now construct in $\AA_S$ an ascending sequence of subobjects of $X$.  Let $A_0 = S \otimes \Hom(S,X)$ and denote $C_0 = \coker(S \otimes\Hom(S,X) \longrightarrow X)$.  We will assume $C_0 \not\cong 0$.

We choose a decomposition $C_0 \cong X_1 \oplus C'_0$ where $X_1$ is indecomposable, hence by Lemma \ref{lemma:Cone} we know $S$ occurs both as subobject and as quotient object of $X_1$ and of every direct summand of $C'_0$.  Consider the following diagram with exact rows and columns
$$\xymatrix{
&& 0\ar[d] & 0\ar[d] & \\
0 \ar[r] & A_0 \ar[r]\ar@{=}[d] & A_1 \ar[r]\ar[d] & S \otimes \Hom(S,X_1) \ar[r]\ar[d] & 0 \\
0 \ar[r] & A_0 \ar[r] & X \ar[r]\ar[d] & X_1 \oplus C'_0 \ar[r]\ar[d] & 0 \\
&& C_1 \oplus C'_0\ar@{=}[r]\ar[d] & C_1 \oplus C'_0 \ar[d] & \\
&& 0 & 0 &
}$$

It follows from Lemma \ref{lemma:Cone} that $S$ occurs both as subobject and quotient of every indecomposable of $C_1 \oplus C'_0$ where $C_1 = \coker (S \otimes\Hom(S,X_1) \longrightarrow X_1)$.  Hence, using $C_0 \not\cong 0$, we have found a subobject $A_1 \in \AA_S$ of $X$ such that $l_S(A_0) < l_S(A_1)$.  Iteration and using that the length is bounded, we see that $X \in \AA_S$.
\end{proof}

\begin{remark}
It follows from previous proposition that all Auslander-Reiten components of $\AA$ are \emph{homogeneous tubes}, i.e. they are of the form $\bZ A_\infty / \langle \tau \rangle$, cfr. Figure \ref{figure:Tube}, were the bottom element is endo-simple.
\end{remark}

Finally, we will formulate a useful corollary.

\begin{corollary}\label{corollary:TubesDirected}
Every cycle $X_0 \to X_1 \to \cdots \to X_n \to X_0$ of non-zero non-isomorphisms between indecomposable objects belongs to a single homogeneous tube.
\end{corollary}

\begin{proof}
Directly from Propositions \ref{proposition:EndoSimplesDirected} and \ref{proposition:Building}.
\end{proof}

\begin{remark}
It follows that the set of homogeneous tubes of the category $\AA$ are directed, thus there can be no cycle containing two objects from different homogeneous tubes.
\end{remark}
\section{Classification}

Let $\AA$ be a connected $k$-linear abelian Ext-finite 1-Calabi-Yau category.  In this section, we wish to classify all such categories up to derived equivalence.  If every two endo-simples of $\AA$ are isomorphic, then $\AA$ is equivalent to the finite dimensional nilpotent representations of the one loop quiver.

So, assume there are at least two non-isomorphic endo-simples, $E$ and $B$.  By connectedness and Proposition \ref{proposition:Building}, we assume $\Hom(E,B) \not= 0$. First, we will find a $t$-structure in $\Db \AA$ such that the heart $\HH$ admits an ample sequence $\EE$.  Then we will use Theorem \ref{theorem:Polishchuk} to show $\AA \cong \qgr R(\EE)$.  A discussion of $R(\EE)$ will then complete the classification of abelian 1-Calabi-Yau categories up to derived equivalence.

From here on, we will always denote $\Hom(E,B)$ by $V$ and its dimension by $d$.

\subsection{The sequence $\EE$ and a $t$-structure in $\Db \AA$}\label{subsection:t}

With $E$ and $B$ as above, associate the autoequivalence $t=T_B : \Db \AA \longrightarrow \Db \AA$ and the sequence $\EE = (E_i)$ where $E_i = t^i E$.

The following will define a $t$-structure in $\CC$ with a hereditary heart $\HH$.
\begin{eqnarray*}
\ind D^{\leq 0} &=& \{X \in \ind \CC \mid \mbox{there is a path from $E_i$ to $X$, for an $i \in \bZ$}\} \\
\ind D^{\geq 1} &=& \ind \CC \setminus \ind D^{\leq 0}
\end{eqnarray*}

If follows directly from this definition that $t$ restricts to an autoequivalence on $\HH$, which we will also denote by $t$.  Note that this implies $E_i \in \Ob \HH$, for all $i \in \bZ$.  Also, since $\Hom(B[-1],E_i) \not= 0$, there is no path from $E_i$ to $B[-1]$ and hence we have $B \in \Ob \HH$.

It follows from Theorem \ref{theorem:Berg} that $\HH$ is hereditary and $\Db \HH \cong \Db \AA$.  Since $\HH$ is a 1-Calabi-Yau category, the results we have proved about $\AA$ apply to $\HH$ as well.

Note that, since $t^i B \cong B$, we find there is a natural isomorphism $\Hom(E,B) \cong \Hom(E_i,B)$ and as such, we get triangles of the form $B[-1] \otimes V^\ast \longrightarrow E_{i-1} \longrightarrow E_i \longrightarrow B \otimes V^\ast$.  Such a triangle in $\Db \AA$ gives rise to an exact sequence
$$0 \longrightarrow E_{i-1} \longrightarrow E_i \longrightarrow B \otimes V^\ast \longrightarrow 0$$
in $\HH$, which is the universal extension of $E_{i-1}$ with $B$ and all these exact sequences lie in the same $t$-orbit.

Using Proposition \ref{proposition:EndoSimplesDirected}, we may prove following easy lemma.

\begin{lemma}\label{lemma:E} Let $\EE = (E_i)_{i\in I}$ and $B$ as above, then
\begin{enumerate}
\item $\Hom(E_i,E_j) = \Ext(E_j,E_i) = 0$ for $i>j$,
\item $\Hom(B,E_i) = \Ext(B,E_i) = 0$ for all $i \in I$.
\end{enumerate}

\end{lemma}

If $\HH$ is of the form $\coh X$ for an elliptic curve $X$ (which we will show below to be the case) one may verify that $E$ corresponds to a stable vector bundle of rank $\dim V$ and $B$ to the structure sheaf $k(P)$ of a point $P$. The $E_i$ are equal to $E(-iP)$.

\subsection{$\EE$ is an ample sequence in $\HH$}

We now wish to show the sequence $\EE = (E_i)_{i \in \bZ}$ is ample.  The following lemma will be useful.

\begin{lemma}\label{lemma:EiX}
If $\Hom(E_i,X) \not=0$, then $\Hom(E_j,X) \not=0$ for all $j \leq i$.
\end{lemma}

\begin{proof}
It suffices to show that $\Hom(E_{i-1},X) \not=0$.  Since $\Hom(E_i,X) \not=0$ and $t$ is an auto-equivalence, we know $\Hom(E_{i-1},t^{-1}X) \not=0$.  Applying the functor $\Hom(E_{i-1},-)$ to the exact sequence
$$0 \longrightarrow t^{-1} X \longrightarrow  X \longrightarrow  B \otimes \Hom(X,B)^\ast \longrightarrow 0$$
yields $\Hom(E_{i-1},X) \not= 0$.
\end{proof}

\begin{proposition}\label{proposition:EE}
In $\HH$ the sequence $\EE = (E_i)$ is ample.
\end{proposition}

\begin{proof}
First, we will show $\EE$ is projective.  Therefore, let $X \to Y$ be an epimorphism and let $K$ be the kernel.  By the construction of $\HH$ in \S\ref{subsection:t}, we know there are paths from the sequence $\EE$ to every direct summand of $K$.  Hence, by Corollary \ref{corollary:TubesDirected}, we know $\Hom(K, E_i)=0$ for $i<<0$ and, by the Calabi-Yau property, $\Ext(E_i,K)=0$.  Thus $\Hom(E_i,X) \to \Hom(E_i,Y)$ is surjective for $i < n$.

Next, we will show $\EE$ is coherent.  Thus we consider an object $X \in \HH$ and we may assume there is a $j \in \bZ$ such that $\Hom(E_{j+2},X) \not= 0$, and hence by Lemma \ref{lemma:EiX}, that $\Hom(E_i,X) \not= 0$ for all $i < j+2$.  Fix an $i < j$, we will prove that $f: E_{i-1}\longrightarrow X$ factors through $E_i \oplus E_{j}$.  Iteration then implies $f$ factors through a number of copies of $E_{j-1} \oplus E_j$, and hence $\EE$ is coherent.

To prove previous claim, it will be convenient to work in the derived category. The following two triangles in $\Db \HH$ will be used
\begin{equation}\label{equation:Triangle1}
\xymatrix@1{B \otimes V^\ast [-1] \ar[r]^-{\theta} & E_{i-1} \ar[r] & E_i \ar[r] & B \otimes V^\ast}
\end{equation}
and
\begin{equation}\label{equation:Triangle2}
\xymatrix@1{B \otimes V^\ast [-1] \ar[r]^-{\varphi} & E_{j} \ar[r] & E_{j+1} \ar[r] & B \otimes V^\ast}
\end{equation}
where $V = \Hom(E_i,B) \cong \Hom(E_{j+1},B)$.  We may assume $f: E_{i-1}\longrightarrow X$ does not factor though $E_i$, hence from triangle (\ref{equation:Triangle1}) it follows that the composition $f \circ \theta \not= 0$.

Note that, since $\Hom(E_{j+1},X) \not= 0$, we may use Corollary \ref{corollary:TubesDirected} to see $\Hom(X,E_{j+1}) = 0$, and hence also $\Ext(E_{j+1},X)=0$.

Applying the functor $\Hom(-,X)$ on triangle (\ref{equation:Triangle2}) and using $\Ext(E_{j+1},X)=0$, shows that every map $B \otimes V^\ast [-1] \longrightarrow X$ factors though $\varphi$.  Hence there is a morphism $g : E_{j} \longrightarrow X$ such that the following diagram commutes.
$$\xymatrix{
B \otimes V^\ast [-1] \ar[r]^-{\theta} \ar[d]_{\varphi}& E_{i-1} \ar[d]^f\\
E_{j} \ar[r]^-{g} & X
}$$

Furthermore, applying $\Hom(-,E_{j})$ to triangle (\ref{equation:Triangle1}) yields that $\varphi$ factors through $\theta$, hence there is a map $h : E_{i-1} \longrightarrow E_{j}$ such that $g \circ h \circ \theta = f \circ \theta$, or $(g \circ h - f) \circ \theta = 0$.

Summarizing, $f = g \circ h + f'$, where $f' : E_{i-1} \longrightarrow X$ lies in $\Ker(\theta,X)$ and as such factors through $E_i$.  The map $f$ factors though $E_i \oplus E_{j}$ and we may conclude the sequence $\EE$ is coherent.

Finally, we show the sequence $\EE$ is ample.  Let $X$ be an indecomposable object.  Due to the construction of $\HH$, we know that there is an oriented path from $E_n$ to $X$, for a certain $n \in \bZ$.  Thus it suffices to prove that if $\Hom(E_n,X) \not= 0$, then there is a finite set $I \subset \bZ$ such that 
$$\bigoplus_{i\in I} E_i \otimes \Hom(E_i,X) \longrightarrow X$$
is an epimorphism.

Let $i_1, \ldots, i_m \in \bZ$ be as in the definition of coherence.  Consider the map
\begin{equation}\label{equation:Ample}
\theta : \bigoplus_{j=1}^{m} E_{i_j} \otimes \Hom(E_{i_j},X) \longrightarrow X
\end{equation}
and let $C = \coker \theta$.  To ease notation, we will refer to the domain of $\theta$ by $\dom \theta$.

There is an exact sequence $0 \longrightarrow \im \theta \longrightarrow X \longrightarrow C \longrightarrow 0$.  Using the Calabi-Yau property, we see $\Hom(\im \theta, C) \not= 0$, and since $\im \theta$ is a quotient object of $\dom \theta$, this yields $\Hom(\dom \theta, C) \not= 0$.  Hence we may assume there is an $i_j$ such that $\Hom(E_{i_j}, C) \not= 0$.

Since $\EE$ is projective, there is an $l<<0$ such that the induced map in $\Hom(E_l,C)$ lifts to a map in $\Hom(E_l,X)$.  Again using coherence, this map should factor through $\dom \theta$.  We may conclude $C=0$, and hence $\theta$ is an epimorphism.
\end{proof}

\subsection{Description of $R = R(\EE)$}

Having shown in Proposition \ref{proposition:EE} that $\EE$ is an ample sequence, we may invoke Proposition \ref{theorem:Polishchuk} to see the that $\HH \cong \cohproj A(\EE)$.

We will now proceed to discuss the graded algebra $R = R(\EE)$.  In particular, we wish to show $R$ is a finitely generated domain of Gelfand-Kirillov dimension 2 which admits a Veronese subalgebra generated in degree one.  It would then follow from \cite{Artin95} that $R$ is noetherian and that $\qgr R$ is equivalent to $\coh X$ where $X$ is a curve, while it would follow from Corollary \ref{corollary:Polishchuk} that $\HH \cong \qgr R$.

We start by showing $\GKdim R = 2$.

\begin{lemma}
Let $\EE= (E_i)_{i \in I}$ and $B$ be as before.  If $j > i$, then
$$\dim \Hom(E_i,E_j) = (j-i)d^2$$
where $d=\dim\Hom(E_0,B)$.
\end{lemma}

\begin{proof}
We apply $\Hom(E_i,-)$ to the short exact sequence
$$0 \longrightarrow E_{j-1} \longrightarrow E_{j} \longrightarrow B \otimes \Hom(E_0,B)^* \longrightarrow 0.$$
We will proceed by induction on $j > i$.  Note that $\dim \Hom(E_i,B) = \dim \Hom(E_0,B)^* = d$ and Lemma \ref{lemma:E} implies that $\Ext(E_i,E_j)=0$.

If $j=i+1$, then it follows from $\dim \Hom(E_i,E_i) = \dim\Ext(E_i,E_i)=1$ that $\dim \Hom(E_i,E_j) = d^2$.
For higher $j$, we find by induction $\dim \Hom(E_i,E_j) = (j-i)d^2$.
\end{proof}


\begin{lemma}\label{lemma:Domain}
Assume $E$ and $B$ are non-isomorphic endo-simple objects of $\Db \AA$ chosen such that $d=\dim \Hom_{\Db \AA}(E,B)$ is minimal and $d \not= 0$.  Then $R$ is a domain.
\end{lemma}

\begin{proof}
It suffices to show every non-zero non-isomorphism $f : E_0 \longrightarrow E_i$ is a monomorphism.  We will prove this by induction on $i$.  The case $i=0$ is trivial.  So let $i \geq 1$.

Since $\im f$ is a quotient object of $E_0$ and $\dim \Hom(E,B)=d$, we see that $\dim \Hom(\im f, B) \leq d$, and due to the minimality of $d$, we know that either $\dim \Hom(\im f, B)=0$, or $\dim \Hom(\im f, B)=d$ and $\im f$ is an endo-simple object.

If $\dim \Hom(\im f, B)=0$, the inclusion $\im f \hookrightarrow E_i$ has to factor through a map $j:\im f \longrightarrow E_{i-1}$.

$$\xymatrix{
&&E_0 \ar@{->>}[d]\\
&&{\im f} \ar@{^{(}->}[d]\ar[ld]_{\exist j}\\
0\ar[r]&E_{i-1}\ar[r]&E_i\ar[r]&B\otimes\Hom(E_i,B)^*\ar[r]&0}$$  

Composition gives a non-zero map $E_0 \longrightarrow E_{i-1}$ which is a monomorphism by the induction hypothesis.  We conclude that $f$ is a monomorphism.

We are left with $\dim \Hom(\im f,B)=d$, and hence $\dim \Hom(K,B)=0$ where $K = \ker f$.  With $\EE$ being ample, we may assume there is a $k \in \bZ$, such that $E_k$ maps non-zero to every direct summand of $K$.  Using the exact sequence $0 \to E_k \to E_{k+1} \to B \otimes \Hom(E_{k+1},B)^* \to 0$, we find that for every $l \in \bZ$, $E_l$ maps non-zero to every direct summand of $K$.  Hence $\Hom(K,E_i) = 0$ and thus $K=0$.  We conclude that $f$ is a monomorphism.
\end{proof}

In general, however, $R$ will not be generated in degree 1.  We show that the Veronese subalgebra $R^{(3)} = \oplus_k R_{3k}$ of $R$ is generated in degree 1.

\begin{lemma}\label{lemma:Degree1}
The sequence $\EE^{(3)}=(E_{3k})_{k\in \bZ}$ is an ample sequence.  Furthermore $R^{(3)} = R(\EE^{(3)})$ is generated in degree 1.
\end{lemma}

\begin{proof}
The sequence $\EE^{(3)}$ is projective and ample since $\EE$ is.  Coherence of $\EE^{(3)}$ may be shown as in the proof of Proposition \ref{proposition:EE}.

Next, we prove $R^{(3)}$ is generated in degree one.  Therefore, it suffices to show that for every $k>1$ every map $E_0 \to E_{3k}$ factors through the canonical map $\theta : E_0 \to E_3 \otimes \Hom(E_0,E_3)^\ast$.  We write $V = \Hom(E_0,E_3)$ and we consider the triangle
$$\xymatrix@1{C \ar[r] & {E_0} \ar[r]^-{\theta} & E_3 \otimes V^\ast \ar[r]& C[1]}$$
where $C = T_{E_3} E_0$ is an endo-simple object since $T_{E_3}$ is an automorphism.  Applying the functor $\Hom(-,E_{3k})$ to this triangle gives the exact sequence
$$0 \longrightarrow \Hom(C[1],E_{3k}) \longrightarrow \Hom(E_3 \otimes V^\ast,E_{3k}) \longrightarrow \Hom(E_0,E_{3k}) \longrightarrow \Hom(C,E_{3k}) \longrightarrow 0.$$
We now consider the dimensions of these vector spaces.  Since 
$$\dim \Hom(E_0,E_{3k}) = (3k)d^2 < \dim \Hom(E_3 \otimes V^\ast,E_{3k}) = 9(k-1)d^4$$
we may see $\Hom(C[1],E_{3k}) \not= 0$ and $\dim \Hom(C,E_{3k}) \not= \dim \Hom(C[1],E_{3k})$, hence $E_{3k} \not\cong C[1]$.

Using Proposition \ref{proposition:EndoSimplesDirected}, we obtain $\Hom(C,E_{3k})=0$, hence every map $E_0 \longrightarrow E_{3k}$ lifts through $\theta$ and the algebra $R^{(3)}$ is generated in degree one.
\end{proof}

\subsection{Classification up to derived equivalence}
We are now ready to prove the main result of this article.

\begin{theorem}\label{theorem:Main}
Let $\AA$ be a connected $k$-linear abelian Ext-finite 1-Calabi-Yau category, then $\AA$ is derived equivalent to either
\begin{enumerate}
\item the category of finite dimensional representations of $k[[t]]$, or 
\item the category of coherent sheaves on an elliptic curve $X$.
\end{enumerate}
\end{theorem}

\begin{proof}
By Proposition \ref{proposition:EndoSimple} we know there are endo-simple objects.  First, assume all endo-simple objects are isomorphic.  Using Proposition \ref{proposition:Building} we easily see that $\AA$ is equivalent to $\Modfd k[[t]]$.

Next, assume there are at least two non-isomorphic endo-simple objects.  Since $\AA$ is connected and using Proposition \ref{proposition:Building}, we may choose two endo-simples, $E$ and $B$, such that $\Hom(E,B) \not= 0$, yet with a minimal dimension.  Let $\HH$ be the abelian category constructed in \S\ref{subsection:t}.

By Lemmas \ref{lemma:Domain} and \ref{lemma:Degree1}, we know $R^{(3)} = R(\EE^{(3)})$ is a domain of GK-dimension 2 which is finitely generated by elements of degree one, hence by \cite{Artin95} we find that $R^{(3)}$ is noetherian and $\qgr R^{(3)}$ is equivalent to the coherent sheaves on a curve $X$.

Since $R$ is noetherian, it follows from \ref{theorem:Polishchuk} that $\HH$ is equivalent to $\qgr R^{(3)}$.

The structure sheaf $\OO_X$ of $X$ is an endo-simple object.  Since the genus of $X$ is $\dim H^1(\OO_X) = \dim \Ext(\OO_X, \OO_X) = \dim \Hom(\OO_X,\OO_X) = 1$, we know $X$ is an elliptic curve.
\end{proof}

\subsection{Classification of abelian categories}\label{subsection:Abelian}

We will now combine Theorem \ref{theorem:Main} with \cite[Proposition 5.1]{Burban06} to obtain a description of all abelian 1-Calabi-Yau categories.  First, we recall some results from \cite{Happel96}.

Let $\AA$ be any hereditary abelian category.  A \emph{torsion theory} on $\AA$, $(\FF,\TT)$, is a pair of full additive subcategories of $\AA$, such that $\Hom(\TT,\FF)=0$ and having the additional property that for every $X \in \Ob \AA$ there is a short exact sequence
$$0 \longrightarrow T \longrightarrow X \longrightarrow F \longrightarrow 0$$
with $F \in \FF$ and $T \in \TT$.

We will say the torsion theory $(\FF,\TT)$ is \emph{split} if $\Ext(\FF,\TT)=0$.  In case of a split torsion theory we obtain, by \emph{tilting}, a hereditary category $\HH$ derived equivalent to $\AA$ with an induced split torsion theory $(\TT,\FF[1])$.  Furthermore, the category $\HH$ will only be hereditary if and only if $(\FF,\TT)$ is a split torsion theory.

We now discuss all possible torsion theories when $\AA$ is equivalent to $\coh X$.  Note that, since $\HH$ will be 1-Calabi-Yau and hence hereditary, all torsion theories on $\AA$ will be split.

Let $(\FF,\TT)$ be a torsion theory on $\AA$, and let $\EE$ be an indecomposable of $\TT$.  Then every indecomposable $\FF$ with slope strictly larger than $\mu(\EE)$ has to be in $\TT$ since $\Hom(\EE,\FF) \not= 0$.  Furthermore, if $\EE$ is in $\TT$ and there is a path from $\EE$ to an indecomposable $\EE'$, then $\EE' \in \ind \TT$. 

We may now give a characterization of all possible torsion theories.

\begin{theorem}\label{theorem:Abelian}\cite{Burban06}
Let $X$ be an elliptic curve.  Every category $\HH$ derived equivalent to $\AA = \coh X$ may be obtained by tilting with respect to a torsion theory.  Moreover, all torsion theories on $\coh X$ are split and may be described as follows.  Let $\theta \in \bR \cup \{\infty\}$.  Denote by $\AA_{> \theta}$ and $\AA_{\geq \theta}$ the subcategory of $\AA$ generated by all indecomposables $\EE$ with $\mu(\EE) > \theta$ and $\mu(\EE) \geq \theta$, respectively.  All full subcategories $\TT$ of $\AA$ with $\AA_{\geq \theta} \subseteq \TT \subseteq \AA_{> \theta} \subseteq \AA$ give rise to a torsion theory $(\FF,\TT)$, with $\ind \FF = \ind \AA \setminus \ind \TT$.
\end{theorem}

\begin{proof}
That these are all possible torsion theories, follows from the above discussion.  That all categories $\HH$ may be obtained in this way, is shown in \cite[Proposition 5.1]{Burban06}.  Alternatively, it is straightforward to check these torsion theories generate all bounded $t$-structures on $\Db \AA$ up to shifts.
\end{proof}

\begin{example}
We give some examples of torsion theories.  In here $\HH$ always stands for the category tilted with respect to the described torsion theory.
\begin{enumerate}
\item If $\theta \in \bQ \cup \{\infty\}$ and $\TT = \AA_{> \theta}$, then the tilted category $\HH$ is equivalent to $\coh X$.  If $\TT = \AA_{\geq \theta}$, then $\HH$ is dual to $\AA$.
\item If $\theta \in \bR \setminus \bQ$ and $\TT = \AA_{> \theta} = \AA_{\geq \theta}$ then $\HH$ is equivalent to the category of holomorphic bundles on a noncommutative two-torus (\cite{Polishchuk04}).
\end{enumerate}
\end{example}

Theorem \ref{theorem:Abelian} classifies all categories derived equivalent to $\coh X$.  We further need to classify all categories derived equivalent to $\BB = \Modfd k[[t]]$.

Let $\HH$ be such a category derived equivalent to $\BB$.  Then $\HH$ induces a $t$-structure $(D^{\geq 0}, D^{\leq 0})$ on $\Db \BB$.  Since this $t$-structure is split, we may assume the heart $\HH = D^{\leq 0} \cap D^{\geq 0}$ contains the endo-simple object $E$ of $\BB[0]$ and, since $\BB$ has only one endo-simple object, this is the unique endo-simple object of $\HH$, up to isomorphism.

Moreover, for every $X \in \BB$ we have $\Hom(X,B) \not=0$ and $\Hom(B,X) \not=0$, thus we have $\BB \subseteq D^{\leq 0} \cap D^{\geq 0} = \HH$.

Since $\BB$ has only one endo-simple object, $E$ is the unique endo-simple object of $\HH$, up to isomorphism.  From this we infer $\BB = \HH$ as subcategories of $\Db \BB$.

We conclude that every category derived equivalent to $\Modfd k[[t]]$ is in fact equivalent to $\Modfd k[[t]]$.

\bibliographystyle{amsplain}

\providecommand{\bysame}{\leavevmode\hbox to3em{\hrulefill}\thinspace}
\providecommand{\MR}{\relax\ifhmode\unskip\space\fi MR }
\providecommand{\MRhref}[2]{%
  \href{http://www.ams.org/mathscinet-getitem?mr=#1}{#2}
}
\providecommand{\href}[2]{#2}

%

\begin{thebibliography}{10}

\bibitem{Artin95}
M.~Artin and J.~T. Stafford, \emph{Noncommutative graded domains with quadratic
  growth}, Invent. Math. \textbf{122} (1995), no.~2, 231--276.

\bibitem{Atiyah57}
M.~F. Atiyah, \emph{Vector bundles over an elliptic curve}, Proc. London Math.
  Soc. (3) \textbf{7} (1957), 414--452.

\bibitem{Beilinson82}
A.~A. Be{\u\i}linson, J.~Bernstein, and P.~Deligne, \emph{Faisceaux pervers},
  Analysis and topology on singular spaces, I (Luminy, 1981), Ast\'erisque,
  vol. 100, Soc. Math. France, Paris, 1982, pp.~5--171.

\bibitem{Berg}
Carl~Fredrik Berg and Adam-Christiaan van Roosmalen, \emph{Projective
  components in hereditary abelian categories satisfying {S}erre duality}, in
  preperation.

\bibitem{Burban06}
Igor Burban and Bernd Kreussler, \emph{Derived categories of irreducible
  projective curves of arithmetic genus one}, Compositio Math. \textbf{142}
  (2006), 1231--1262.

\bibitem{Happel01}
Dieter Happel, \emph{A characterization of hereditary categories with tilting
  object}, Invent. Math. \textbf{144} (2001), no.~2, 381--398.

\bibitem{Happel96}
Dieter Happel, Idun Reiten, and Smal{\o}~Sverre O., \emph{Tilting in abelian
  categories and quasitilted algebras}, Mem. Amer. Math. Soc. \textbf{120}
  (1996), no.~575, viii+ 88.

\bibitem{Keller06}
Bernhard Keller and Idun Reiten, \emph{Acyclic {C}alabi-{Y}au categories},
  preprint (2006).

\bibitem{Polishchuk04}
A.~Polishchuk, \emph{Classification of holomorphic vector bundles on
  noncommutative two-tori}, Doc. Math. \textbf{9} (2004), 163--181.

\bibitem{Polishchuk05}
\bysame, \emph{Noncommutative proj and coherent algebras}, Math. Res. Lett.
  \textbf{12} (2005), no.~1, 63--74.

\bibitem{Polishchuk03}
A.~Polishchuk and A.~Schwarz, \emph{Categories of holomorphic vector bundles on
  noncommutative two-tori}, Comm. Math. Phys. \textbf{236} (2003), no.~1,
  135--159.

\bibitem{ReVdB02}
I.~Reiten and M.~Van~den Bergh, \emph{Noetherian hereditary abelian categories
  satisfying {S}erre duality}, J. Amer. Math. Soc. \textbf{15} (2002), no.~2,
  295--366.

\bibitem{Ringel05}
Claus~Michael Ringel, \emph{Hereditary triangulated categories}, Compos. Math.,
  to appear.

\bibitem{Seidel01}
Paul Seidel and Richard Thomas, \emph{Braid group actions on derived categories
  of coherent sheaves}, Duke Math. J. \textbf{108} (2001), no.~1, 37--108.

\bibitem{vanRoosmalen06}
Adam-Christiaan van Roosmalen, \emph{Classification of abelian hereditary
  directed categories satisfying {S}erre duality}, Trans. Amer. Math. Soc.
  \textbf{360} (2008), no.~5, 2467--2503.

\end{thebibliography}

\end{document}